\documentclass[12pt,a4paper]{amsart}
\usepackage[latin1]{inputenc}
\usepackage{latexsym}
\usepackage{color,graphicx,shortvrb}
\usepackage{amsmath, amssymb}
\usepackage{amsfonts}
\usepackage[colorlinks, bookmarks=true]{hyperref}

\newtheorem{theorem}{Theorem}[section]
\newtheorem{lemma}[theorem]{Lemma}
\newtheorem{proposition}[theorem]{Proposition}
\newtheorem{corollary}[theorem]{Corollary}

\theoremstyle{definition}
\newtheorem{definition}[theorem]{Definition}
\newtheorem{remark}[theorem]{Remark}

\author{J. M. Almira$^*$, Kh. F. Abu-Helaiel}
\title{On solutions of $f(x)+f(a_1x)+\cdots+f(a_Nx)=0$ and related equations}
\thanks{$^*$ Corresponding author}

\begin{document}


\subjclass[2010]{}


\begin{abstract} We study some properties of the solutions of the functional equation $$f(x)+f(a_1x)+\cdots+f(a_Nx)=0,$$ which was introduced in the literature by  Mora, Cherruault and Ziadi in 1999, for the case $a_k=k+1$, $k=1,2,\cdots,N$ \cite{mora_1} and studied by Mora \cite{mora_2} and Mora and Sepulcre \cite{mora_3, mora_4}.   \end{abstract}

\maketitle

\markboth{J. M. Almira, Kh. F. Abu-Helaiel}{On solutions of $f(x)+f(a_1x)+\cdots+f(a_Nx)=0$}

\section{Motivation}
The functional equation  \begin{equation}\label{mora}
 f(x)+f(a_1x)+\cdots+f(a_Nx)=0,\end{equation}
was introduced in the literature by  Mora, Cherruault and Ziadi in 1999, for the case $a_k=k+1$, $k=1,2,\cdots,N-1$ \cite{mora_1} and studied by Mora \cite{mora_2} and Mora and Sepulcre \cite{mora_3,mora_4}. Concretely, the equations $f(x)+f(2x)=0$ and $f(x)+f(2x)+f(3x)=0$ were used for modeling certain processes related to combustion of hydrogen in a car engine \cite{mora_1, mora_2} and, later, the most general equation
\begin{equation} \label{ecumora}
f(x)+f(2x)+\cdots+f(Nx)=0
\end{equation}
was studied from a more theoretical point of view \cite{mora_2,mora_3,mora_4}. Concretely, by imposing a solution of the form $z^\alpha$, these authors have shown that there exists a strong connection between the continuous solutions of $(\ref{ecumora})$ and the zeroes of the exponential functions $G_N(z)=1+2^z+\cdots+N^z$, and they have developed a very interesting theory with some deep results concerning the distribution of the zeroes of $G_N(z)$. Note that the functions $H_N(z)=G_N(-z)$ are the partial sums of the Riemann zeta function $\zeta(z)=\sum_{n=1}^\infty n^{-z}$, so that the study of the zeroes of $G_N(z)$ is an important problem, connected with the well known Riemann's conjecture.

In \cite{mora_2} the author proved that  if $f(z)$ is a solution of $(\ref{ecumora})$ and  $f(z)\neq 0$, then $f(z)$ cannot be analytic at $z=0$. Furthermore, he also proved that the set of solutions of $(\ref{ecumora})$ that are analytic on $\Omega=\mathbb{C}\setminus (-\infty,0]$ is an infinite dimensional vector space. In this paper we study the more general equation  $(\ref{mora})$. In Section 2 we prove that if $f(x)$  is a solution of  $(\ref{mora})$, with $0<a_1<\cdots<a_N$ and $a_k\neq 1$ for $k=1,2,\cdots,N$, there exists a positive natural number $m=m(a_1,\cdots,a_N)$ such that, for any $\delta >0$, if $f\in \mathbf{C}^{(m)}[0,\delta]$, then $f_{|[0,\delta]}=0$ and, if if $f\in \mathbf{C}^{(m)}[-\delta,0]$, then $f_{|[-\delta,0]}=0$. In particular, if
$f\in \mathbf{C}^{(k)}(\mathbb{R})$ is a solution of $(\ref{mora})$ and $k\geq m(a_1,\cdots,a_N)$, then $f=0$. We also give upper and lower bounds for $m(a_1,\cdots,a_N)$. In Section 3 we concentrate on the study of the equation  $(\ref{mora})$ with the additional restriction $x>0$. We transform this equation  into the easier one
\begin{equation}   \label{nueva}
g(w)+g(w+b_1)+\cdots+g(w+b_N)=0 \ \ (w\in\mathbb{R}),
\end{equation}
and we give a new elegant argument to prove that the set of continuous solutions of both equations is an infinite dimensional vector space. Furthermore, we also study the existence of continuous periodic solutions for  $(\ref{nueva})$.

\section{A regularity result for $f(x)+f(a_1x)+\cdots+f(a_Nx)=0$}
We study, for $N\geq 1$, the functional equation:
$$
 f(x)+f(a_1x)+\cdots+f(a_Nx)=0,
$$
where $0<a_1<a_2<\cdots<a_N$ are positive real numbers, $a_k\neq 1$ for all $k\in\{1,\cdots,N\}$, and $x$ is a real variable. Moreover, in all what follows we set $a_0=1$.

The first thing we observe is that  we can assume that $1<a_1$ since otherwise, making the change of variable $y=a_1x$, the equation $(\ref{mora})$ is transformed into the equation
\begin{equation*} 
f(y)+f(b_1y)+\cdots+ \cdots +f(b_Ny)=0 \ \ (x\in (0,\infty)),
\end{equation*}
where
\[
(b_1,b_2,\cdots,b_N)==\left\{
\begin{array}{llll}
(\frac{a_2}{a_1},\cdots, \frac{a_k}{a_1},\frac{1}{a_1},\frac{a_{k+1}}{a_1},\cdots ,\frac{a_N}{a_1}) & \  & \text{if} &   a_k< 1< a_{k+1} \\
(\frac{a_2}{a_1},\cdots, \frac{a_N}{a_1},\frac{1}{a_1}) & \  & \text{if} & a_N < 1 \\
\end{array}
\right.
.
\]
Hence, we impose the conditions $0<a_0=1<a_1<\cdots<a_N$ through this paper. Furthermore, we use the notation $\mathbf{a}=(a_1,a_2,\cdots,a_N)$.

\begin{definition} Given $I\subseteq \mathbb{R}$ an interval, we say that a function $f:I\to\mathbb{R}$ satisfies $(\ref{mora})$ (or that $f$ is a solution of $(\ref{mora})$ on $I$) if  $f(x)+f(a_1x)+\cdots+f(a_Nx)=0$ whenever $\{x,a_1x,\cdots,a_Nx\}\subseteq I$.
\end{definition}

 \begin{definition} Let $N\geq 1$ and $\mathbf{a}=(a_1,a_2,\cdots,a_N)$. We define the natural number
 \[
 m(\mathbf{a})=\min\left\{m\in\mathbb{N}: \sum_{k=0}^{N-1}(\frac{a_k}{a_N})^{m}<1\right\}.
 \]
 \end{definition}

Obviously, $m(\mathbf{a})$ is well defined because all fractions $\frac{a_k}{a_N}$ appearing under the summation symbol satisfy $0<\frac{a_k}{a_N}<1$, and the number $N$ of summands is fixed.

\begin{theorem}\label{teoregularidad}  Let $\delta>0$ be a positive real number. Then:
\begin{itemize}
\item[$(a)$] If $f$ is a solution of  $(\ref{mora})$ on $[0,\delta]$ and $f\in \mathbf{C}^{(m(\mathbf{a}))}[0,\delta]$, then $f_{|[0,\delta]}=0$.
\item[$(b)$] If $f$ is a solution of  $(\ref{mora})$ on $[-\delta,0]$ and  $f\in \mathbf{C}^{(m(\mathbf{a}))}[-\delta,0]$, then $f_{|[-\delta,0]}=0$.
\end{itemize}
\end{theorem}

\noindent \textbf{Proof. } We only prove part $(a)$ of the theorem, since part $(b)$ follows with the very same arguments. Assume that $f\in \mathbf{C}^{(m(\mathbf{a}))}[0,\delta]$  is a solution of $(\ref{mora})$ over $[0,\delta]$. Let us set $x=h/a_N$, so that equation $(\ref{mora})$ is transformed into
  \begin{equation}\label{mora_transf}
 f(h/a_N)+f(a_1h/a_N)+\cdots+f(a_{N-1}h/a_N)+f(h)=0.\end{equation}
 Taking derivatives $m(\mathbf{a})$ times at $(\ref{mora_transf})$ and defining $\varphi(h)=f^{(m(\mathbf{a}))}(h)$, we get
  \begin{equation*}
 \varphi(h)=(-1)\left[(\frac{1}{a_N})^{m(\mathbf{a})}\varphi(\frac{h}{a_N})+ (\frac{a_1}{a_N})^{m(\mathbf{a})}\varphi(\frac{a_1}{a_N}h)+\cdots+ (\frac{a_{N-1}}{a_N})^{m(\mathbf{a})}\varphi(\frac{a_{N-1}}{a_N}h)\right].\end{equation*}
 It follows that
 \[
 \|\varphi\|_{[0,\delta]}\leq \|\varphi\|_{[0,\delta]} \sum_{k=0}^{N-1}(\frac{a_k}{a_N})^{m(\mathbf{a})},
 \]
so that $\|\varphi\|_{[0,\delta]}=0$ (since  $\sum_{k=0}^{N-1}(\frac{a_k}{a_N})^{m(\mathbf{a})}<1$). This implies that $f_{|[0,\delta]}$ is a polynomial of degree $\leq m(\mathbf{a})-1$.  In particular, $f$ is real analytic on $[0,\delta]$ and, given  any natural number  $m\in \mathbb{N}$, we have that
 \begin{equation*}
 f^{(m)}(h)=(-1)\left[(\frac{1}{a_N})^{m}f^{(m)}(\frac{h}{a_N})+ (\frac{a_1}{a_N})^{m}f^{(m)}(\frac{a_1}{a_N}h)+\cdots+ (\frac{a_{N-1}}{a_N})^{m}f^{(m)}(\frac{a_{N-1}}{a_N}h)\right], \end{equation*}
 so that
  \begin{equation*}
 f^{(m)}(0)\left[1+ (\frac{1}{a_N})^{m}+ (\frac{a_1}{a_N})^{m}+\cdots+ (\frac{a_{N-1}}{a_N})^{m}\right]=0.\end{equation*}
 Hence $f^{(m)}(0)=0$ for all $m$. This means that $f_{|[0,\delta]}=0$. {\hfill $\Box$}

\begin{corollary} If $f\in \mathbf{C}^{(k)}(\mathbb{R})$ is a solution of $(\ref{mora})$ and $k\geq m(\mathbf{a})$, then $f=0$. In particular, $f=0$ is the unique solution of $(\ref{mora})$ which admits infinitely many derivatives in all the real line.
\end{corollary}

\begin{remark}It is important to note that, for all $N\geq 2$, the function $f:[0,\infty)\to\mathbb{R}$ given by $f(x)=0$ can be extended in infinitely many ways to a solution $\widetilde{f}$ of equation $(\ref{ecumora})$ such that  $\widetilde{f}\in \mathbf{C}(\mathbb{R})\cap \mathbf{C}^{(\infty)}(\mathbb{R}\setminus \{0\})$ and $\widetilde{f}$ does not vanish identically.
To prove this, we use that, for $N\geq 2$, the function $G_N(z)=1+2^z+\cdots+N^z$ has infinitely many zeros in the complex plane  (see \cite[Proposition 1]{mora_2} for a proof of this fact). Thus, let $\alpha=a+\mathbf{i}b\in\mathbb{C}$ be a zero of $G_N(z)$ and let us define $\widetilde{f}_{\alpha}(x)=\mathbf{Re}(|x|^{\alpha})$ for $x<0$ and $\widetilde{f}(x)=0$ for $x\geq 0$. Then $\widetilde{f}_{\alpha}$ is clearly a solution of equation $(\ref{ecumora})$ in $[0,\infty)$, and, for $x<0$ we have that
\[
\sum_{k=1}^N\widetilde{f}_{\alpha}(kx)=\mathbf{Re}\left(|x|^{\alpha}\sum_{k=1}^N k^{\alpha}\right) = 0,
\]
so that $\widetilde{f}_{\alpha}$ is also a solution of equation $(\ref{ecumora})$ in $\mathbb{R}$. \end{remark}

\begin{remark} The natural number  $m(\mathbf{a})$ appearing in Theorem \ref{teoregularidad}  is not optimal. For example, it is clear that, for $N=1$, $a_1=2$,  $m(\mathbf{a})=1$. Now, let us assume that $f$ is a solution of $f(x)+f(2x)=0$ and $f(h_0)\neq 0$ for some $h_0>0$. Then $f(\frac{1}{2}h_0)=(-1)f(h_0)$ and, for all $k\in\mathbb{N}$ we have that $f(\frac{1}{2^k}h_0)=(-1)^kf(h_0)$. This obviously implies that $f(x)$ cannot be continuous at $x=0$.
\end{remark}

The following result gives a simple estimation of the number $m(\mathbf{a})$:

 \begin{proposition}  $$\frac{1}{2}\frac{a_N}{\max_{0\leq k<N}(a_{k+1}-a_k)}-1\leq m(\mathbf{a})\leq \frac{a_N}{\min_{0\leq k<N}(a_{k+1}-a_k)}$$
 \end{proposition}

 \noindent \textbf{Proof. } The result follows directly from the interpretation of the sum $$\sum_{k=0}^{N-1}\left(\frac{a_k}{a_N}\right)^m\left(\frac{a_{k+1}}{a_N}-\frac{a_k}{a_N}\right)$$
 as a lower Riemann sum for the integral $\int_0^1x^mdx$ and the sum $$\sum_{k=0}^{N}\left(\frac{a_k}{a_N}\right)^m\left(\frac{a_{k}}{a_N}-\frac{a_{k-1}}{a_N}\right)$$
 as an upper Riemann sum for the same integral (here we impose $a_{-1}=0$). Concretely, we have that
 \[
 \sum_{k=0}^{N-1}\left(\frac{a_k}{a_N}\right)^m\left(\frac{a_{k+1}}{a_N}-\frac{a_k}{a_N}\right)\leq \int_0^1x^mdx=\frac{1}{m+1},
 \]
 so that, for all $m\geq \frac{a_N}{\min_{0\leq k<N}(a_{k+1}-a_k)}$,
 \begin{eqnarray*}
  \sum_{k=0}^{N-1}\left(\frac{a_k}{a_N}\right)^m &\leq &  \sum_{k=0}^{N-1}\left(\frac{a_k}{a_N}\right)^m\left(\frac{a_{k+1}}{a_N}-\frac{a_k}{a_N}\right) \frac{a_N}{\min_{0\leq k<N}(a_{k+1}-a_k)}\\
  & \leq & \frac{1}{m+1}\frac{a_N}{\min_{0\leq k<N}(a_{k+1}-a_k)}<1,
 \end{eqnarray*}
 which implies that $m(\mathbf{a})\leq \frac{a_N}{\min_{0\leq k<N}(a_{k+1}-a_k)}$.

 On the other hand,
 \begin{eqnarray*}
 \sum_{k=0}^{N}\left(\frac{a_k}{a_N}\right)^m &\geq & \sum_{k=0}^{N}\left(\frac{a_k}{a_N}\right)^m\left(\frac{a_{k}}{a_N}-\frac{a_{k-1}}{a_N}\right)\frac{a_N}{\max_{0\leq k<N}(a_{k+1}-a_k)}\\
 & \geq &  (\int_0^1x^mdx) \frac{a_N}{\max_{0\leq k<N}(a_{k+1}-a_k)} \\
 &=& \frac{1}{m+1} \frac{a_N}{\max_{0\leq k<N}(a_{k+1}-a_k)},
 \end{eqnarray*}
 so that, for all $m < \frac{1}{2}\frac{a_N}{\max_{0\leq k<N}(a_{k+1}-a_k)}-1$,
 \[
  \sum_{k=0}^{N-1}\left(\frac{a_k}{a_N}\right)^m  \geq \frac{1}{m+1}\frac{a_N}{\max_{0\leq k<N}(a_{k+1}-a_k)}-1>1,
  \]
 which implies that $m(\mathbf{a})\geq \frac{1}{2}\frac{a_N}{\max_{0\leq k<N}(a_{k+1}-a_k)}-1$.  {\hfill $\Box$}

 Note that the best bounds for $m(\mathbf{a})$ appear when the points $a_k$ are equidistributed (i.e., when $a_k=1+k(a_1-1)$ for all $k\in\{0,1,\cdots,N\}$) since it is precisely in this case when  $\max_{0\leq k<N}(a_{k+1}-a_k)$ attains its minimum and $\min_{0\leq k<N}(a_{k+1}-a_k)$ attains its maximum, and both coincide with $d=a_1-1$. Hence, in this case we get the following bounds: $$\frac{1}{2}N-1\leq m((1+d,1+2d,\cdots,1+Nd))\leq N,$$ which are independent of the separation $d$.

\section{A related functional equation}

Let us consider the equation
\begin{equation}\label{morageneral}
f(x)+f(a_1x)+\cdots +f(a_Nx)=0 \ \ (x\in (0,\infty)),
\end{equation}
where $a_0=1<a_1<a_2<\cdots<a_N$ are real numbers.
If we set $x=e^w$ and $g(w)=f(e^w)$ then, taking into account that $f(a_kx)=f(a_ke^w)=f(e^{w+\ln a_k})=g(w+\ln a_k)$, the equation $(\ref{morageneral})$ can be written as
\begin{equation}\label{aditiva}
g(w)+g(w+b_1)+\cdots+g(w+b_N)=0 \ \ (w\in\mathbb{R}),
\end{equation}
where $0<b_k=\ln a_k <\ln a_{k+1}< b_{k+1}$, $k=1,\cdots,N-1$.

\begin{lemma} \label{auxiliar} Let us assume that $g:[0,b_N]\to \mathbb{R}$ is a continuous function which satisfies $(\ref{aditiva})$. Then there exists a unique $\widetilde{g}\in\mathbf{C}(\mathbb{R})$ such that $\widetilde{g}$ is a solution of  $(\ref{aditiva})$ on $\mathbb{R}$ and $\widetilde{g}_{|[0,b_N]}=g$.
\end{lemma}
\noindent \textbf{Proof. } First of all, we note that $g:[0,b_N]\to \mathbb{R}$ satisfies $(\ref{aditiva})$ if and only if it satisfies the interpolation condition:
\[
g(0)+g(b_1)+\cdots+g(b_N)=0.
\]
Moreover, if $y,w$ denote two real numbers satisfying the relation  $y=w+b_N$ and $\widetilde{g}$ is any solution of  $(\ref{aditiva})$ on $\mathbb{R}$, then a simple substitution shows that $\widetilde{g}$ satisfies
\[
\widetilde{g}(y)+\widetilde{g}(y-b_N)+\widetilde{g}(y-(b_N-b_1))+\cdots+\widetilde{g}(y-(b_N-b_{N-1}))=0 \ \ (y\in\mathbb{R}).
\]
We will use this relation to (uniquely) define on $[0, \infty)$ the solution $\widetilde{g}$ such that $\widetilde{g}_{|[0,b_N]}=g$. Furthermore, we will use the original equation $(\ref{aditiva})$  to (uniquely) extend the solution $\widetilde{g}$ over the negative part of the real axis.

Let us set $I_0=[0,b_N]$, $\widetilde{g}_0=g$, and define, for $h\in I_1=[b_N, b_N+(b_{N}-b_{N-1})]$, the function
\[
\widetilde{g}_1(y)= (-1)\left[\widetilde{g}_0(y-b_N)+\widetilde{g}_0(y-(b_N-b_1))+\widetilde{g}_0(y-(b_N-b_2))+\cdots+\widetilde{g}_0(y-(b_N-b_{N-1}))\right].
\]
Obviously, $\widetilde{g}_1$ is well defined, since $t\in I_1$ implies that
$$0\leq y-b_N\leq y-(b_N-b_{1})\leq y-(b_N-b_{2})\leq \cdots \leq y-(b_N-b_{N-1})\leq b_N.$$
Moreover, $\widetilde{g}_1\in \mathbf{C}(I_1)$.
For $k\geq 2$, we set  $I_{k}=[b_N+(k-1)(b_{N}-b_{N-1}), b_N+k(b_{N}-b_{N-1})]$ and
\[
\widetilde{g}_k(y)= (-1)\left[\widetilde{g}_{k-1}(y-b_N)+\widetilde{g}_{k-1}(y-(b_N-b_1))+\cdots+\widetilde{g}_{k-1}(y-(b_N-b_{N-1}))\right] \ \ (y\in I_{k}).
\]
Let us now consider the negative part of the real axis. Set $I_{-1}=[-b_1,0]$ and
\[
\widetilde{g}_{-1}(x)= (-1)\left[\widetilde{g}_{0}(x+b_1)+\widetilde{g}_{0}(x+b_2)+\cdots+\widetilde{g}_{0}(x+b_N)\right] \ \ (x\in I_{-1}).
\]
For $k\leq -2$, we set  $I_{k}=[kb_{1}, (k+1)b_1]$ and
\[
\widetilde{g}_{k}(x)= (-1)\left[\widetilde{g}_{k+1}(x+b_1)+\widetilde{g}_{k+1}(x+b_2)+\cdots+\widetilde{g}_{k+1}(x+b_N)\right] \ \ (x\in I_{k}).
\]
Clearly, $\bigcup_{k\in\mathbb{Z}} I_k=\mathbb{R}$ and $\widetilde{g}(x)=\widetilde{g}_k(x)$ ($x\in I_k$, $k\in\mathbb{Z}$) is the function we were looking for. The uniqueness is guaranteed by the construction we have used for the definition of $\widetilde{g}$.

{\hfill $\Box$}

For the proof of the following theorem, we need firstly to recall the concept of exponential polynomial which is of common use for people working on functional equations.

\begin{definition} We say that $f(x)\in\mathbf{C}(\mathbb{R})$ is a (real) exponential polynomial if $f(x)$ is the solution of some ordinary homogeneous linear differential equation with constant coefficients, $y^{(n)}+a_1y^{(n-1)}+\cdots+a_ny=0$. These functions are completely characterized as finite $\mathbb{R}$-linear combinations of the real and imaginary parts of functions of the form $m(x)=x^k e^{\lambda x}$, where $k\leq n-1$ is a natural number and $\lambda$ is a complex number. Furthermore, they can also be characterized as the continuous solutions of certain functional equations which do not involve the use of derivatives, such as Popoviciu's equation (see, for example, \cite{roscau}, \cite{rado}):
\[
\det \left[
\begin{array}{cccccc}
f(x) & f(x+h) & \cdots & f(x+nh) \\
f(x+h) & f(x+2h) & \cdots & f(x+(n+1)h)   \\
f(x+2h) & f(x+3h) &  \cdots & f(x+(n+2)h)  \\
\vdots &  & \ddots & \vdots &  \\
f(x+nh) & f(x+(n+1)h) & \cdots & f(x+2nh)
\end{array}
\right] = 0.
\]
We say that $f(x)\in\mathbf{C}(\mathbb{R},\mathbb{C}):=\{f:\mathbb{R}\to\mathbb{C},\  f\text{ is continuous}\}$ is a (complex) exponential polynomial if $f(x)$ is a finite $\mathbb{C}$-linear combination of functions of the  form $m(x)=x^k e^{\lambda x}$, where $k$ is a natural number and $\lambda$ is a complex number.
\end{definition}

\begin{theorem} \label{dimension} Let $\mathbf{S}=\{g\in\mathbf{C}(\mathbb{R}):g \text{ is a solution of } (\ref{aditiva})\}$. Then $\mathbf{S}$ is an infinite dimensional vector space.  As a consequence, the space of continuous solutions of $(\ref{morageneral})$  is also an infinite dimensional vector space.
\end{theorem}

\noindent \textbf{Proof. }  It is well known (see \cite{anselone}, \cite{engert}) that, if $V$ is a finite dimensional subspace of $\mathbf{C}(\mathbb{R},\mathbb{C})$ and $V$ is invariant by translations (i.e., $f(x)\in V$ implies $g_L(x)=f(x-L)\in V$ for all $L\in\mathbb{R}$), then all elements of $V$ are (complex) exponential polynomials. In particular, all elements of $V$ are analytic functions. On the other hand, the space $\mathbf{S}$ is obviously invariant by translations and can be considered in a natural way as a subspace of $\mathbf{C}(\mathbb{R},\mathbb{C})$. Moreover, Lemma \ref{auxiliar} implies that $\mathbf{S}$ is nonempty. Hence, the proof will end as soon as we find a continuous solution of  (\ref{aditiva}) which is not an exponential polynomial.

Let us define
$$
g(x)= \left\{
\begin{array}{lll}
1 & \text{ if }   0\leq  x \leq b_{N-1} \\
\frac{b_N+Nb_{N-1}-(N+1)x}{b_N-b_{N-1}} &   \text{ if }   b_{N-1}<x\leq b_N \end{array}
\right.
$$
Obviously, $g$ satisfies $(\ref{aditiva})$ in $[0,b_N]$, so that Lemma \ref{auxiliar} implies that there exists $\widetilde{g}\in\mathbf{S}$ such that $\widetilde{g}_{|[0,b_N]}=g$. In particular, $\widetilde{g}$ is not an exponential polynomial, since $g$ is not differentiable.  {\hfill $\Box$}

It is interesting to note that, in some cases, the function $\widetilde{g}$ we get in the construction shown at the proof of Theorem  \ref{dimension} is periodic. For example, if we impose $b_k=k$, $k=1,2,\cdots,N$ and we follow all steps of the proof, we get that $\widetilde{g}(x)$ is the $(N+1)$-periodic extension of the function :
$$
g(x)= \left\{
\begin{array}{lll}
1 & \text{ if }   0\leq  x \leq N-1 \\
-(N+1)x + N^2  &   \text{ if }   N-1<x\leq N \\
(N+1)x -N(N+2) &   \text{ if }   N<x\leq N+1 \\
\end{array}
\right.
$$

Thus, an interesting question is, under which conditions on $(b_1,\cdots,b_N)$ can we guarantee that equation $(\ref{aditiva})$ admits periodic solutions?   The following theorem partially solves this question:


\begin{theorem} The equation $(\ref{aditiva})$ admits a continuous periodic solution $g\neq 0$ if and only if there exists $\alpha\in\mathbb{R}$ such that
\begin{equation}\label{condiciongeneral}
\left\{ \begin{array}{lll}
1+\sum_{k=1}^N\cos(\alpha b_k) & = & 0 \\
\sum_{k=1}^N\sin(\alpha b_k) &  =& 0   \\
\end{array}
\right.
\end{equation}
Furthermore, in such a case, there are trigonometric polynomials satisfying equation  $(\ref{aditiva})$ which are periodic of fundamental period equal to $2\pi/\alpha$. Finally, the equation $(\ref{aditiva})$ admits continuous periodic solutions for $(b_1,\cdots,,b_N)$ if and only if it admits continuous periodic solutions for $(d b_1,\cdots, d b_N)$ for all $d> 0$.
\end{theorem}
\noindent \textbf{Proof. } Assume $g$ is a continuous $T$-periodic solution of  $(\ref{aditiva})$ and set  $\theta=2\pi/T$. Let
\[
g(x)=\sum_{k=1}^{\infty}(a_k(g)\cos(k\theta t)+b_k(g)\sin(k\theta t)) +\frac{a_0(g)}{2}
\]
be the Fourier series expansion of $g$ (this expansion exists because, being $g$ continuous, its Fourier coefficients are well defined). Then $h(x)=g(x)+\sum_{k=1}^Ng(x+b_k)$ is also continuous and $T$-periodic, so that its Fourier coefficients are well defined. Indeed, a simple computation gives  $a_0(h)=a_0(g)$ and
\[
\left(
\begin{array}{lll}
a_k(h) \\
b_k(h) \\
\end{array}
\right)
= \left(
\begin{array}{lll}
1+\sum_{k=1}^N\cos(k\theta b_k) & \sum_{k=1}^N\sin(k\theta b_k) \\
-\sum_{k=1}^N\sin(k\theta b_k) &   1+\sum_{k=1}^N\cos(k\theta b_k) \\
\end{array}
\right) \left(
\begin{array}{lll}
a_k(g) \\
b_k(g) \\
\end{array}
\right)  \text{ for all } k\geq 1.
\]
Now, the function $h$ vanishes identically if and only if all its Fourier coefficients are zero. Thus, if $g$ is a solution of $(\ref{aditiva})$  then  $a_0(g)=0$ and
\begin{equation} \label{matriz}
\left(
\begin{array}{lll}
1+\sum_{k=1}^N\cos(k\theta b_k) & \sum_{k=1}^N\sin(k\theta b_k) \\
-\sum_{k=1}^N\sin(k\theta b_k) &   1+\sum_{k=1}^N\cos(k\theta b_k) \\
\end{array}
\right)
\left(
\begin{array}{lll}
a_k(g) \\
b_k(g) \\
\end{array}
\right)  = \left(
\begin{array}{lll}
0 \\
0 \\
\end{array}
\right)
\end{equation}
for all $k\geq 1$. Let us denote by $A_k$ the matrix appearing in $(\ref{matriz})$. If $g$ is not the zero function, then there exists $k\geq 1$ such that $(a_k(g),b_k(g))\neq (0,0)$ and, for this concrete value of $k$ we should have
\[
\det(A_k)=\left(1+\sum_{k=1}^N\cos(k\theta b_k) \right)^2 +  \left(\sum_{k=1}^N\sin(k\theta b_k) \right)^2 = 0.
\]
In other words, the system of equations $(\ref{condiciongeneral})$ is satisfied for $\alpha=k\theta$.  What is more: as soon as $\det(A_k)=0$ we have that $A_k=0$ is the null matrix, which implies that for all  $(a_k,b_k)\in\mathbb{R}^2$ the trigonometric polynomial
$$g(x) = a_k\cos(k\theta t)+b_k\sin(k\theta t)$$
satisfies $(\ref{aditiva})$ and is a periodic function with fundamental period $T=\frac{2\pi}{k\theta}$.  The last claim of the theorem follows from the fact that $\alpha$ is a solution of  $(\ref{condiciongeneral})$ for $(b_1,\cdots,b_N)$ if and only if $\frac{\alpha}{d}$ if a solution of $(\ref{condiciongeneral})$ for $(d b_1,\cdots,d b_N)$.

{\hfill $\Box$}

\begin{corollary} Let $d>0$ and set $b_k=k d$, $k=1,2,\cdots N$ and let $m\in\mathbb{Z}\setminus (N+1)\mathbb{Z}$.  Then equation  $(\ref{aditiva})$ admits continuous periodic solutions of period $T=\frac{N+1}{m}$. \end{corollary}

\noindent \textbf{Proof. } It is only necessary to make the proof for $d=1$. Assume that $b_k=k$ for $k=1,2,\cdots,N$. Then $(\ref{condiciongeneral})$ becomes:
\begin{equation} \label{condiciongeneralcaso}
\left\{ \begin{array}{lll}
\frac{\sin\left( \frac{N+1}{2}\alpha \right)}{\sin\left( \frac{\alpha}{2}\right)} \cos\left( \frac{N}{2}\alpha \right)  = & 0 \\
\frac{\sin\left( \frac{N+1}{2}\alpha \right)}{\sin\left( \frac{\alpha}{2}\right)} \sin\left( \frac{N}{2}\alpha \right)    = & 0   \\
\end{array}
\right. .
\end{equation}
Hence, in this case, to find $\alpha\in\mathbb{R}$ which solves  $(\ref{condiciongeneral})$ is equivalent to find a real solution $\alpha$ of
\[
\left(\frac{\sin\left( \frac{N+1}{2}\alpha \right)}{\sin\left( \frac{\alpha}{2}\right)}\right)^2 = \left[\frac{\sin\left( \frac{N+1}{2}\alpha \right)}{\sin\left( \frac{\alpha}{2}\right)} \cos\left( \frac{N}{2}\alpha \right) \right]^2+ \left[\frac{\sin\left( \frac{N+1}{2}\alpha \right)}{\sin\left( \frac{\alpha}{2}\right)} \sin\left( \frac{N}{2}\alpha \right) \right]^2=0.
\]
Equivalently, we are looking for the real solutions of
\[
\frac{\sin\left( \frac{N+1}{2}\alpha \right)}{\sin\left( \frac{\alpha}{2}\right)}=0,
\]
which exist and are given by $\alpha = \frac{2m\pi}{N+1}$, $m\in\mathbb{Z}\setminus (N+1)\mathbb{Z}$. This ends the proof. {\hfill $\Box$}

\begin{corollary} Let us assume that $b> 0$. The equation \begin{equation}\label{dos} g(x)+g(x+a)+g(x+b)=0\ \ (x\in\mathbb{R})\end{equation} admits a continuous periodic solution $g\neq 0$ if and only if
\begin{equation} \label{condiciondos}
\frac{a}{b}  \in\{\frac{2+3k}{1+3m}, \frac{1+3m}{2+3k}: (m,k)\in\mathbb{Z}^2\}.
\end{equation}
In particular, if $a/b\in\mathbb{R}\setminus\mathbb{Q}$ then $(\ref{dos})$ admits no continuous periodic solutions. Furthermore, there are infinitely many rational numbers $p/q$ such that $a/b=p/q$ implies that $(\ref{dos})$ admits no continuous periodic solutions. \end{corollary}

\noindent \textbf{Proof. } In this case, the equations we must study are given by:  \begin{equation}\label{condiciongeneraltres}
\left\{ \begin{array}{lll}
1+ cos(\alpha a) +cos(\alpha b) & = & 0 \\
sin(\alpha a) +sin(\alpha b)  &  =  & 0   \\
\end{array}
\right.
\end{equation}
Solving the second equation in the system we get that $\alpha b= -\alpha a + 2k\pi$ or $\alpha b= \alpha a + (2k+1)\pi$ for a certain $k\in\mathbb{Z}$. We consider both cases separately:

\noindent \textbf{Case 1: $\alpha b= \alpha a + (2k+1)\pi$. } Introducing the corresponding values into the first equation in the system, we get
\[
1+ cos(\alpha a) +cos(\alpha a + (2k+1)\pi) =1\neq 0.
\]
Hence, in this case we get no solutions of $(\ref{condiciongeneraltres})$.

\noindent \textbf{Case 2: $\alpha b= -\alpha a + 2k\pi$. } Introducing the corresponding values into the first equation in the system, we get
\[
 1+ cos(\alpha a) +cos(-\alpha a + 2k\pi) =1 +2 cos(\alpha a) = 0 \Leftrightarrow
\alpha a \in \{2\pi/3+2m\pi, \pi/3+(2m+1)\pi\}_{m\in\mathbb{Z}}.
\]
Taking into account that $\alpha b= -\alpha a + 2k\pi$, we conclude that $\alpha$ is a solution of  the system of equations $(\ref{condiciongeneraltres})$ if and only if it is a solution of at least one of the following two  systems:
\begin{equation}\label{I}
\left\{ \begin{array}{lll}
\alpha a  & = & 2\pi/3+2m\pi \\
\alpha b  &  = & -(2\pi/3+2m\pi) +2k\pi  \\
\end{array}
\right.      \text{ for some } (m,k)\in\mathbb{Z}^2.
\end{equation}
or
\begin{equation}\label{II}
\left\{ \begin{array}{lll}
\alpha a  & = & \pi/3+(2m+1)\pi\\
\alpha b  &  = & -(\pi/3+(2m+1)\pi) +2k\pi  \\
\end{array}
\right.      \text{ for some } (m,k)\in\mathbb{Z}^2.
\end{equation}
Let us study the equations given in $(\ref{I})$. Clearly, we should have $a,b\neq 0$. Furthermore,
\[
\alpha = \frac{1}{a} \left(2\pi/3+2m\pi\right) = \frac{1}{b} \left(-(2\pi/3+2m\pi) +2k\pi\right)
\]
In particular, the system has real solutions if and only if
\[
\frac{b}{a} \in \{\frac{-1+3(k-m)}{1+3m}: (m,k)\in\mathbb{Z}^2\} =\{\frac{2+3k}{1+3m}: (m,k)\in\mathbb{Z}^2\}.
\]
As the parameters $a,b$ are interchangeable in all the argument above, we conclude that
\[
\frac{b}{a} \in \{\frac{2+3k}{1+3m}, \frac{1+3m}{2+3k}: (m,k)\in\mathbb{Z}^2\}.
\]
Finally, it is clear that the system $(\ref{II})$ is a particular case of $(\ref{I})$.

{\hfill $\Box$}

\begin{remark} The condition $(\ref{condiciondos})$ can be studied for any particular instance of $a,b$. For example, if $a=b$ then $(\ref{condiciondos})$ becomes
\[
1  \in\{\frac{2+3m}{1+3k}, \frac{1+3k}{2+3m}: (k,m)\in\mathbb{Z}^2\}
\]
which is clearly impossible. It follows that there are no continuous periodic functions $g\neq 0$ satisfying $g(x)+2g(x+a)=0$.
We give here a direct proof of this fact:  Assume, on the contrary, that $g(x)+2g(x+a)=0$ and $g(x)$ is continuous and $T$-periodic with $T>0$ a fundamental period. Let $s_0\in [0,T]$ be such that $|g(s_0)|\neq 0$. There are two possibilities:
\begin{itemize}
\item[$(a)$] The numbers $\{a,T\}$ are linear dependent when dealing $\mathbb{R}$ as a $\mathbb{Q}$-vector space.
    In this case, there are natural numbers $k,m$ such that $kT=ma$. Then
    \[
    0<|g(s_0)|=|g(s_0+kT)|=|g(s_0+ma)|=\frac{1}{2^m}|g(s_0)|<|g(s_0)|,
    \]
    which is a contradiction.
\item[$(b)$] $\dim \mathbf{span}_{\mathbb{Q}}\{a,T\}=2$. Taking into account that $g$ is continuous and periodic, it is uniformly continuous. Hence, given $\varepsilon>0$ there exists $\delta>0$ such that $|g(x)-g(y)|<\varepsilon$ whenever $|x-y|<\delta$. Now, our hypothesis on $\{a,T\}$ implies that for any $\delta>0$ there exists $n,m\in\mathbb{Z}$ such that $|nT+ma|<\delta$ (indeed, we may assume $m>0$). Hence, if $d=|g(s_0)|>0$ and $\varepsilon<d/2$,
    \[
    |g(s_0+nT+ma)-g(s_0)|<\varepsilon
    \]
    which implies $|g(s_0+nT+ma)|>d/2$. On the  other hand,
    \[
    |g(s_0+nT+ma)|=|g(s_0+ma)|=\frac{1}{2^m}|g(s_0)| \leq \frac{1}{2}d,
    \]
    a contradiction.
\end{itemize}
Of course, a direct proof of existence (or nonexistence)  of continuous periodic solutions of equation $(\ref{dos})$ for each instance of the parameters $a,b$, would be a difficult task. Instead of that, checking if condition $(\ref{condiciondos})$ is (or it is not) satisfied is always easy.
\end{remark}

\section{Acknowledgements}

We are very grateful to the referee for reading our paper carefully
and thoroughly, and making many helpful suggestions.


 \bibliographystyle{amsplain}


\bigskip

\footnotesize{J. M. Almira

Departamento de Matem\'{a}ticas. Universidad de Ja\'{e}n.

E.P.S. Linares,  C/Alfonso X el Sabio, 28

23700 Linares (Ja\'{e}n) Spain

Email: jmalmira@ujaen.es }



\vspace{1cm}

\footnotesize{Khader. F. Abu-Helaiel

Departamento de Estadística e Investigación Operativa. Universidad de Ja\'{e}n.

Campus de Las Lagunillas

23071 Ja\'{e}n, Spain

Email: kabu@ujaen.es
}



\end{document}